\documentclass[a4paper,12pt,final]{amsart}
\usepackage{amssymb}
\usepackage{showkeys}
\usepackage[pdfauthor={Eugen Stumpf},%
pdftitle={A note on local center manifolds for differential equations with state-dependent delay},%
colorlinks=true]{hyperref}
\newcommand{\N}{\mathbb{N}}  % The natural numbers.
\newcommand{\R}{\mathbb{R}}  % The real numbers.
\newcommand{\C}{\mathbb{C}}  % The complex numbers.
\DeclareMathOperator{\id}{id} % The identity map
\newtheorem{thm}{Theorem}

\title[A Note on Local Center Manifolds]{A Note on Local Center
Manifolds for Differential Equations with State-Dependent Delay}

\author{Eugen Stumpf}
\address{Department of Mathematics, University of Hamburg, Bundesstrasse 55, 20146 Hamburg, Germany}
\email{eugen.stumpf@math.uni-hamburg.de}

\urladdr{www.math.uni-hamburg.de/home/stumpf/index\_en.html}



{}
\copyrightinfo{2015}{Eugen Stumpf}

\PII{}

\begin{document}

\begin{abstract}
  In this note we consider local invariant manifolds of functional differential equations
  $x^{\prime}(t)=f(x_{t})$ representing differential equations with state-dependent delay.
  Starting with a local center-stable and a local center-unstable manifold of the functional
  differential equation at a stationary point, we construct, by a straightforward application
  of the Implicit Mapping Theorem, a local center manifold.
\end{abstract}

\maketitle

\section{Introduction}
Although the first studies of differential equations with state-dependent delay go back
at least to the very beginning of 19th century, a significant research interest in that type of equations began not till the second half of the last century. In the light of this, it is unsurprising that only in the last descent there were developed a general framework for studying differential equations with state-dependent delay in the context of the dynamical systems theory. The starting point here was the work \cite{Walther2003} of Walther proving, under some mild smoothness assumptions, the existence of a continuous semiflow with $C^{1}$-smooth time-$t$-maps for differential equations with state-dependent delay.
In subsequent years this semiflow and its properties were studied in different papers under slightly modified smoothness conditions, and the best general reference here is the survey article \cite{Hartung2006} of Hartung et al. and the references therein.

Now it is well known that the mentioned semiflow particularly has different types of local invariant manifolds at a stationary point. For instance, the survey article \cite{Hartung2006} of Hartung et al. contains the construction of so-called \textit{local center manifolds} by employing the Lyapunov-Perron method (with appropriate changes). The $C^{1}$-smoothness of these finite-dimensional invariant manifolds is shown in Krisztin \cite{Krisztin2006}. Some years later the approaches in Hartung et al. \cite{Hartung2006} and in Krisztin \cite{Krisztin2006} were adopted in \cite{Stumpf2011}, in order to prove the existence and $C^{1}$-smoothness of another type of finite-dimensional local invariant manifolds, namely, of so-called \textit{local center-unstable manifolds}. Contrary to that, a proof of $C^{1}$-smooth \textit{local center-stable manifolds}, which are infinite-dimensional, could not be obtained by employing a variation of the Lyapunov-Perron technique. However, as shown in Qesmi and Walther \cite{Qesmi2009}, they arise from local center-stable manifolds of time-$t$-maps of some global semiflows that are modifications of the original semiflow.

In the situation of a flow generated by an ordinary differential equation in a finite-dimensio\-nal space, of course, all three mentioned types of locally invariant manifolds at a stationary point are well understood. In this case, each of them is finite-dimensio\-nal and can be constructed, for instance, by applying the Lyapunov-Perron method as can be found in Vanderbauwhede \cite{Vanderbauwhede1989}, or by using the graph-transformation method as discussed in Kelley \cite{Kelley1967}.
Additionally, it is well known that a local center mani\-fold at a stationary point may also be obtained by intersecting a local center-unstable and a local center-stable manifold at the same stationary point.

Returning to differential equations with state-depen\-dent delay, it suggests itself to ask whether the last point does also hold in the situation of the discussed semiflow; that is, does an intersection between a local center-unstable manifold from \cite{Stumpf2011} and a local center-stable manifold from Qesmi and Walther \cite{Qesmi2009} at the same stationary point form or contain a local center manifold as constructed in Hartung et al. \cite{Hartung2006}? The present paper answers this question in the affirmative and that was to be expected. By using local representations of the involved manifolds as graphs of maps and applying the Implicit Mapping Theorem, we show in a purely analytical way -- in particular, without discussing properties such as transversality etc. from the geometric theory of differentiable manifolds -- that the intersection of a local center-unstable and a local center-stable manifold indeed contains a local center manifold. A simple consequence of this is the fact that given a local center-stable and/or a local center-unstable manifold of the discussed semiflow at a stationary point we always find a local center manifold contained in the considered local center-stable and/or local center-unstable manifold.

The general approach applied here, that is, the construction of a local manifold via intersecting two non-disjoint manifolds by means of the Implicit Mapping Theorem, is certainly not new and it works of course also in the case of a flow induced by an ordinary differential equation. But observe that for a differential equation with state-dependent the situation is somewhat more subtil than for an ordinary differential equation in finite-dimensional space, since the state-space itself is a submanifold with finite codimension of an infinite dimensional Banach space.

The remaining part of this note is organized as follows. The next section contains a brief summary of the general setting. After introducing the differential equation together with the smoothness assumptions under consideration, we discuss here the mentioned semiflow and some of its properties.
In the end, we state the local center manifold theorem -- compare Theorem \ref{thm: CMT} -- obtained in \cite{Hartung2006,Krisztin2006}.

Section 3 forms the main part. Starting with a local center-unstable and a local center-stable manifold, we construct a local center manifold and give so an alternative proof of the local center manifold theorem for differential equations with state-dependent delay.

\section{Outline of a semiflow framework for differential equations with state-dependent delay}
In the sequel we give a short summary of a general concept for studying differential equations with state-dependent delay in the context of the dynamical systems theory. For the left proofs as well as for a thorough discussion of the topic we refer the reader to the survey work \cite{Hartung2006} of Hartung at al. and the references therein.

Let $n\in\N$ and $h>0$ be fixed. Further, let $\|\cdot\|_{\R^{n}}$ stand for any fixed norm on the $n$-dimensional Euclidean space $\R^{n}$, and $C$ for the Banach space of all continuous functions $\varphi:[-h,0]\to\R^{n}$ equipped with the usual norm $\|\varphi\|_{C}:=\sup_{-h\leq s\leq 0}\|\varphi(s)\|_{\R^{n}}$. Similarly, we write $C^{1}$ for the Banach space of all continuously differentiable functions $\varphi:[-h,0]\to\R^{n}$
with the norm given by $\|\varphi\|_{C^{1}}:=\|\varphi\|_{C}+\|\varphi^{\prime}\|_{C}$.

Given a function $x:I\to\R^{n}$ defined on some interval $I\subset\R$ and $t\in\R$ with $[t-h,t]\subset I$, we will write $x_{t}$ for the function $[-h,0]\ni s\mapsto x(t+s)\in\R^{n}$, which is also known as the \textit{segment} of $x$ at $t$.

From now on, we consider the functional differential equation
\begin{equation}\label{eq: FDE}
  x^{\prime}(t)=f(x_{t})
\end{equation}
defined by a map $f:C^{1}\supset U\to\R^{n}$ on some open neighborhood $U$ of $0\in C^{1}$ with $f(0)=0$. Under a \textit{solution} of Eq. \eqref{eq: FDE} we understand either a continuously differentiable function $x:[t_{0}-h,t_{+})\to\R^{n}$, $t_{0}<t_{+}\leq \infty$, such that $x_{t}\in U$ for all $t_{0}\leq t<t_{+}$ and $x$ satisfies Eq. \eqref{eq: FDE} as $t_{0}<t<t_{+}$, or a continuously differentiable function $x:\R\to\R^{n}$ satisfying $x_{t}\in U$ and Eq. \eqref{eq: FDE} for all of $t\in\R$. In this sense, the function $x:\R\ni t\mapsto 0\in\R^{n}$ is clearly a solution of Eq. \eqref{eq: FDE} due to the assumption $f(0)=0$.

Before discussing the assumptions on $f$ which will ensure the existence of further solutions of Eq. \eqref{eq: FDE}, we should mention briefly the connection between differential equations with state-dependent delay and equations of the form \eqref{eq: FDE}. For this purpose, consider the differential equation
\begin{equation}\label{eq: delay_equation}
  x^{\prime}(t)=g(x(t-r(x(t))))
\end{equation}
defined by some function $g:\R\to\R^{n}$ with $g(0)=0$ and involving a state-dependent delay given by another function $r:\R^{n}\to[0,h]$. Introducing
\begin{equation*}
  \hat{f}:C^{1}\ni\varphi\mapsto g(\varphi(-r(\varphi(0))))\in\R^{n},
\end{equation*}
we see that
\begin{equation}\label{eq: sample}
  x^{\prime}(t)=g(x(t-r(x(t))))=g(x(t-r(x_{t}(0))))=g(x_{t}(-r(x_{t}(0)))=\hat{f}(x_{t}),
\end{equation}
that is, the differential equation \eqref{eq: delay_equation} with state-dependent delay takes the more abstract form of Eq. \eqref{eq: FDE}. Consequently, instead of studying the original equation \eqref{eq: delay_equation}, we may just as well study Eq. \eqref{eq: sample}.

The presented transformation also works in many other cases of differential equations with state-dependent delay than the simple one discussed above. However, having in mind that Eq. \eqref{eq: FDE} represents a differential equation with a state-dependent delay, we should particularly impose only those smoothness assumptions on $f$, which are typically fulfilled by differential equations with state-depen-dent delay. With this regard, we follow the further development of the ideas contained in Walther \cite{Walther2003}, and suppose that $f$ satisfies the following conditions:
\begin{itemize}
  \item[(S1)] $f$ is continuously differentiable, and
  \item[(S2)] for each $\varphi\in U$ the derivative $Df(\varphi):C^{1}\to\R^{n}$ extends to a linear map $D_{e}f(\varphi):C\to\R^{n}$ such that the map
        $U\times C\ni(\varphi,\psi)\mapsto D_{e}f(\varphi)\psi\in\R^{n}$
      is continuous.
\end{itemize}
Then the closed and, in view of the zero function, nonempty subset
\begin{equation*}
  X_{f}:=\lbrace \varphi\in U\mid \varphi^{\prime}(0)=f(\varphi)\rbrace
\end{equation*}
of $U\subset C^{1}$ forms a continuously differentiable submanifold of $U$ with codimension $n$. For each $\varphi\in X_{f}$, Eq. \eqref{eq: FDE} has a unique (in the forward time-direction) noncontinuable solution $x^{\varphi}:[-h,t_{+}(\varphi))\to\R^{n}$ satisfying $x^{\varphi}_{0}=\varphi$. Furthermore, for all $\varphi\in X_{f}$ and all $0\leq t<t_{+}(\varphi)$ the segments $x^{\varphi}_{t}$ belong to $X_{f}$, and the relations
\begin{equation*}
F(t,\varphi):=x^{\varphi}_{t}
\end{equation*}
as $\varphi\in X_{f}$ and $0\leq t<t_{+}(\varphi)$ define a domain $\Omega\subset[0,\infty)\times X_{f}$ and a continuous semiflow $F:\Omega\to X_{f}$ with $C^{1}$-smooth time-$t$-maps $F_{t}:=F(t,\cdot)$.

Now, recall that $x(t)=0$, $t\in\R$, is a solution of Eq. \eqref{eq: FDE}. For this reason, $\varphi_{0}=0\in X_{f}$ is a stationary point of the semiflow $F$, that is, $F(t,\varphi_{0})=0$ for all $t\geq 0$. The linearization of $F$ at $\varphi_{0}$ is the strongly continuous semigroup $T=\lbrace T(t)\rbrace_{t\geq 0}$ of bounded linear operators $T(t):=DF_{t}(0): T_{0}X_{f}\to T_{0}X_{f}$ on the tangent space
\begin{equation*}
  T_{0}X_{f}:=\left\lbrace \chi\in C^{1}\mid \chi^{\prime}(0)=Df(0)\chi\right\rbrace
\end{equation*}
of $F$ at $\varphi_{0}$, which forms a Banach space with respect to the norm $\|\cdot\|_{C^{1}}$ of $C^{1}$. Given $t\geq 0$ and $\chi\in T_{0}X_{f}$, the operator $T(t)$ maps $\chi$ to the segment $v^{\chi}_{t}$ of the uniquely determined solution $v^{\chi}:[-h,\infty)\to\R^{n}$ of the linear variational equation
\begin{equation*}
  v^{\prime}(t)=Df(0)v_{t}
\end{equation*}
with initial value $v_{0}^{\chi}=\chi$. The infinitesimal generator of $T$ is the linear operator $G:\mathcal{D}(G)\ni \chi\mapsto \chi^{\prime}\in T_{0}X_{f}$
defined on the closed subset
\begin{equation*}
  \mathcal{D}(G):=\left\lbrace\chi\in C^{2}\mid \chi^{\prime}(0)=Df(0)\chi, \chi^{\prime\prime}(0)=Df(0)\chi^{\prime}\right\rbrace
\end{equation*}
of the space $C^{2}$ of all twice continuously differentiable functions from $[-h,0]$ into $\R^{n}$.

In order to describe the spectrum of the linearization $T$, or more precisely, the spectrum $\sigma(G)\subset\C$ of its generator $G$, remember that due to condition (S2) the derivative $Df(0)$ is extendible to a bounded linear operator $D_{e}f(0):C\to\R^{n}$. In particular, $D_{e}f(0)$ defines the linear retarded functional differential equation
\begin{equation*}
  v^{\prime}(t)=D_{e}f(0)v_{t}
\end{equation*}
on $C$. As shown, for instance, in Diekmann et al. \cite{Diekmann1995}, the solutions of the associated Cauchy problems induce a strongly continuous semigroup $T_{e}:=\lbrace T_{e}(t)\rbrace_{t\geq 0}$ of linear bounded operators $T_{e}(t):C\to C$, and $G_{e}:\mathcal{D}(G_{e})\ni \varphi\mapsto \varphi^{\prime}\in C$
with
\begin{equation*}
  \mathcal{D}(G_{e}):=\left\lbrace \chi\in C^{1}\mid \chi^{\prime}(0)=Df_{e}(0)\chi\right\rbrace
\end{equation*}
forms the generator of $T_{e}$. Obviously, $T_{0}X_{f}=\mathcal{D}(G_{e})$. But even more is true: For all $t\geq 0$ and all $\varphi\in \mathcal{D}(G_{e})$ we have $T(t)\varphi=T_{e}(t)\varphi$, and the spectra $\sigma(G)$ and $\sigma(G_{e})$ coincide. The spectrum $\sigma(G_{e})\subset \C$ of $G_{e}$ is given by the zeros of a familiar characteristic function. In particular, it is discrete, contains only eigenvalues with finite dimensional generalized eigenspaces, and for each $\beta\in\R$ the intersection $\sigma(G_{e})\cap \lbrace \lambda\in\C\mid \Re(\lambda)\geq \beta\rbrace$ is either finite or empty. For this reason, each of the spectral parts $\sigma_{u}(G):=\lbrace\lambda\in\sigma(G_{e})\mid\Re(\lambda)>0\rbrace$ and
$\sigma_{c}(G):=\lbrace\lambda\in\sigma(G_{0})\mid \Re(\lambda)=0\rbrace$ is either finite or empty as well. Moreover, it follows that the associated realified generalized eigenspaces $C_{u}$ and $C_{c}$, which are called the \textit{unstable} and the \textit{center space} of $G_{e}$, respectively, are finite dimensional subspaces of $T_{0}X_{f}\subset C^{1}$. In contrast to those, the \textit{stable space} of $G_{e}$, that is, the realified generalized eigenspace given by the eigenvalues with negative real part, is infinite dimensional and belongs not to $T_{0}X_{f}$. In all, these subspaces provide the decomposition
\begin{equation*}
  C=C_{u}\oplus C_{c}\oplus C_{s}
\end{equation*}
of $C$. Since $C^{1}_{s}:=C^{1}\cap C_{s}$ is closed in $C^{1}$ we also have the decomposition
\begin{equation}\label{eq: decomposition of C^1}
  C^{1}=C_{u}\oplus C_{c}\oplus C_{s}^{1}
\end{equation}
of the smaller Banach space $C^{1}$. The unstable and center space of $G$ coincide with $C_{u}$ and $C_{c}$, respectively, whereas the stable space of $G$ is given by the intersection $C_{s}\cap T_{0}X_{f}=C^{1}_{s}\cap T_{0}X_{f}$. In particular, we get
\begin{equation*}
  T_{0}X_{f}=C_{u}\oplus C_{c}\oplus (C_{s}^{1}\cap T_{0}X_{f}).
\end{equation*}

Next, we repeat some facts about local invariant manifolds of $F$ at $\varphi_{0}=0$. For doing so, recall that a \emph{trajectory} of $F$ is a map $\gamma:I\to X_{f}$, $I\subset\R$ an interval, with $\gamma(t)=F(t-s,\gamma(s))$ whenever $s,t\in I$ and $t\geq s$. We begin with a statement contained in Qesmi and Walther \cite[Theorem 1.1]{Qesmi2009}:
There exist an open neighborhood $N_{cs}$ of $0$ in $U$ and a continuously differentiable submanifold $W_{cs}\subset X_{f}$ with
\begin{equation*}
T_{0}W_{cs}=C_{c}\oplus (C^{1}_{s}\cap T_{0}X_{f})
\end{equation*}
which has the following properties.
\begin{itemize}
  \item[(CS 1)] $W_{cs}$ is positively invariant with respect to $F$ relative to $N_{cs}$; that is, for all $\varphi\in W_{cs}$ and all $0\leq t< t_{\varphi}$ with $F([0,t]\times\lbrace\varphi\rbrace)\subset N_{cs}$ we have $F(s,\varphi)\in W_{cs}$ as $0\leq s\leq t$.
  \item[(CS 2)] $W_{cs}$ contains all initial values $\varphi \in X_{f}$ with $F(t,\varphi)\in N_{cs}$ for all $0\leq t<\infty$.
  \item[(CS 3)] If $\gamma:[t,0]\to X_{f}$, $t<0$, is a trajectory of $F$ with $\gamma(0)\in W_{cs}$ and with $\gamma([t,0])\subset N_{cs}$, then $\gamma(s)\in W_{cs}$ as $t\leq s\leq 0$.
\end{itemize}
The submanifold $W_{cs}$ of $X_{f}$ is called a \emph{local center-stable manifold} of $F$ at $\varphi_{0}$. In terms of solutions of Eq. \eqref{eq: FDE}, the assertion (CS 1) means that all the segments of a solution of Eq. \eqref{eq: FDE} with initial value in $W_{cs}$ remain in $W_{cs}$ as long as the solution does not leave the neighborhood $N_{cs}$ of $\varphi_{0}$, whereas assertion (CS 2) says that $W_{cs}$ contains the segments of all sufficiently small solutions $x:[-h,\infty)\to\R^{n}$ of Eq. \eqref{eq: FDE}.

The counterpart of a local center-stable manifold $W_{cs}$ is formed by a so-called \emph{local center-unstable manifold} $W_{cu}$ which does exist under the additional assumption that $\lbrace \lambda\in\sigma(G_{e})\mid\Re(\lambda)\geq 0\rbrace\not=\emptyset$ and so $C_{cu}:=C_{u}\oplus C_{c}\not=\lbrace 0\rbrace$ holds as discussed in \cite[Theorems 1 \& 2]{Stumpf2011}: There exist an open neighborhood $N_{cu}$ of $0$ in $U$ and a continuously differentiable submanifold $W_{cu}\subset X_{f}$ of dimension $\dim C_{cu}\geq 1$ with
\begin{equation*}
  T_{0}W_{cu}=C_{cu}=C_{u}\oplus C_{c}
\end{equation*}
which has the following properties:
\begin{itemize}
  \item[(CU 1)] $W_{cu}$ is positively invariant with respect to $F$ relative to $N_{cu}$; that is, for all $\varphi\in W_{cu}$ and all $0\leq t< t_{\varphi}$ with $F([0,t]\times\lbrace\varphi\rbrace)\subset N_{cu}$ we have $F(s,\varphi)\in W_{cu}$ as $0\leq s\leq t$.
  \item[(CU 2)] If $\gamma:(-\infty,0]\to X_{f}$ is a trajectory of $F$ with $\gamma((-\infty,0])\subset N_{cu}$ then $\gamma(t)\in W_{cu}$ for all $t\leq 0$.
\end{itemize}
Concerning solutions of Eq. \eqref{eq: FDE}, the property (CU 1) of $W_{cu}$ is completely similar to the one of the local center-stable manifold $W_{cs}$, whereas assertion (CU 2) means that $W_{cu}$ contains the segments of all sufficiently small solutions $x:(-\infty,0]\to\R^{n}$ of Eq. \eqref{eq: FDE}. But let us now state the Local Center Manifold Theorem which is proven in \cite{Hartung2006,Krisztin2006} and whose alternative construction is the main ingredient of this note.

\begin{thm}[see Theorem 4.1.1 in \cite {Hartung2006} and Theorem 2.1 in \cite{Krisztin2006}]\label{thm: CMT}
  Suppose that, in addition to the stated assumptions on $f$, $\sigma(G_{e})\cap i\R\not=\emptyset$, that is, $\dim C_{c}\geq 1$ holds. Then there exist open neighborhoods $C_{c,0}$ of $0$ in $C_{c}$ and $C^{1}_{su,0}$ in $C_{s}^{1}\oplus C_{u}$ with $N_{c}=C_{c,0}+C_{su,0}^{1}\subset U$ and a continuously differentiable map $w_{c}:C_{c,0}\to C_{su,0}^{1}$ such that $w_{c}(0)=0$, $Dw_{c}(0)=0$, and such that for the graph
  \begin{equation*}
    W_{c}=\left\lbrace\varphi+w_{c}(\varphi)\mid \varphi\in C_{c,0}\right\rbrace
  \end{equation*}
  of $w_{c}$, which is called a local center manifold of $F$ at $\varphi_{0}=0$, the following holds.
  \begin{itemize}
    \item[(C 1)] $W_{c}\subset X_{f}$, and $W_{c}$ is a $\dim C_{c}$-dimensional continuously differentiable submanifold of $X_{f}$.
    \item[(C 2)] $W_{c}$ is positively invariant with respect to $F$ relative to $N_{c}$; that is, for all $\varphi\in W_{c}$ and all $0\leq t<t_{\varphi}$ with $F([0,t]\times\lbrace\varphi\rbrace)\subset N_{c}$ we have $F(s,\varphi)\in W_{c}$ as $0\leq s\leq t$.
    \item[(C 3)] If $\gamma:\R\to X_{f}$ is a trajectory of $F$ with $\gamma(\R)\subset N_{c}$, then $\gamma(t)\in W_{c}$ for all $t\in\R$.
  \end{itemize}
\end{thm}
In contrast to the discussion of $W_{cs}$ and $W_{cu}$, the above result is formulated in terms of a local representation.  Observe that we have
\begin{equation*}
  T_{0}W_{c}=C_{c}
\end{equation*}
and that by (C 3) the local center manifold $W_{c}$ contains the segments of all globally defined and small enough solutions of Eq. \eqref{eq: FDE}. As already mentioned in the introduction, the original proof of Theorem \ref{thm: CMT} is based on the Lyapunov-Perron method. In the next section, we deduce the statement of Theorem \ref{thm: CMT} from the existence of a local center-stable $W_{cs}$ and a local center-unstable manifold $W_{cu}$ by means of the Implicit Mapping Theorem.

\section{Alternative proof of Theorem \ref{thm: CMT}}
From now on, let the assumptions of Theorem \ref{thm: CMT} be satisfied. Then our discussion in the last section implies that we find both a local center-stable $W_{cs}$ as well as a local center-unstable manifold $W_{cu}$ of $F$ at $\varphi_{0}=0$. Let $N_{cs}$ and $N_{cu}$ denote the corresponding open neighborhoods of $0$ in $U$ where $W_{cs}$ and $W_{cu}$ are positively invariant, respectively. The remaining proof will be divided into five parts as follows: In the first one, we prepare the application of the Implicit Mapping Theorem by introducing a (local) manifold chart for $X_{f}$ and representing both $W_{cs}$ and $W_{cu}$ locally at $\varphi_{0}=0$ as graphs of appropriate maps. The second part then contains the application of the Implicit Mapping Theorem, whereas in the third one we define a local center manifold $W_{c}$. In the final two steps we give a representation of $W_{c}$ as a graph of a continuously differentiable map and show that $W_{c}$ has properties (C 1) -- (C 3) of Theorem \ref{eq: FDE}.

\subsection{\texorpdfstring{Local representations of $X_{f}$, $W_{cs}$ and $W_{cu}$}{Local representations}}
Set $Y:=T_{0}X_{f}$ and fix some subspace $Z\subset C^{1}$ with $\dim Z=n$ such that
\begin{equation*}
  C^{1}=Y\oplus Z
\end{equation*}
holds. Let $P:C^{1}\to C^{1}$ denote the continuous projection of $C^{1}$ along $Z$ onto $Y$. Then recall the decomposition \eqref{eq: decomposition of C^1} of $C^{1}$ and that we have $C_{u},C_{c}\subset Y$. Hence, introducing the closed subspace $Y_{s}:=C^{1}_{s}\cap Y=C_{s}\cap Y$ of $Y$, we obtain
\begin{equation*}
  C^{1}_{s}=Y_{s}\oplus Z
\end{equation*}
and so the additional decomposition
\begin{equation*}
C^{1}=C_{u}\oplus C_{c}\oplus Y_{s}\oplus Z
\end{equation*}
of $C^{1}$. In particular, there are continuous projections $P_{u},P_{c},P_{Y_{s}}:C^{1}\to C^{1}$ of $C^{1}$ onto $C_{u}$, $C_{c}$, and $Y_{s}$, respectively.

Now, observe that we find an open neighborhood $U_{0}\subset U$ of $\varphi_{0}=0$ in $C^{1}$ and an open neighborhood $Y_{0}$ of $0$ in $Y$ such that the equation
\begin{equation*}
  K(\varphi)=P(\varphi-\varphi_{0})=P\varphi
\end{equation*}
defines a manifold chart for $X_{f}$ with $K(\varphi_{0})=0\in Y$ and $K(U_{0}\cap X_{f})=Y_{0}$. The inverse of $K$ is given by a $C^{1}$-smooth map $R:Y_{0}\to U_{0}\subset C^{1}$, and both derivatives $DK(\varphi_{0})$ and $DR(0)$ are equal to the identity operator on $Y$.

Next, we represent both $W_{cs}$ as well as $W_{cu}$ locally at $\varphi_{0}=0$ as a graph of a map defined on some open neighborhood of $0$ in the corresponding tangent space. Consider the submanifold $W_{cs}$ first: There are open neighborhoods $C^{cs}_{cs,0}$ of $0$ in $C_{c}\oplus Y_{s}=T_{0}W_{cs}$ and $C^{cs}_{uz,0}$ of $0$ in $C_{u}\oplus Z$ with
\begin{equation*}
\tilde{N}_{cs}:=C_{cs,0}^{cs}+C_{uz,0}^{cs}\subset N_{cs}\cap U_{0}
\end{equation*}
 and a continuously differentiable map $w_{cs}:C_{cs,0}^{cs}\to C_{uz,0}^{cs}$ with $w_{cs}(0)=0$ and $Dw_{cs}(0)=0$ such that we have
\begin{equation*}
  W_{cs}\cap \tilde{N}_{cs}=\lbrace \varphi+w_{cs}(\varphi)\mid \varphi\in C_{cs,0}^{cs}\rbrace.
\end{equation*}
Similarly, we find open neighborhoods $C^{cu}_{cu,0}$ of $0$ in $C_{c}\oplus C_{u}=T_{0}W_{cu}$ and $C_{sz,0}^{cu}$ of $0$ in $Y_{s}\oplus Z$ with
\begin{equation*}
\tilde{N}_{cu}:=C_{cu,0}^{cu}+ C_{sz,0}^{cu}\subset N_{cu}\cap U_{0}
\end{equation*}
 and a continuously differentiable map $w_{cu}:C_{cu,0}^{cu}\to C_{sz,0}^{cu}$ with $w_{cu}(0)=0$ and $Dw_{cu}(0)=0$ such that
\begin{equation*}
  W_{cu}\cap \tilde{N}_{cu}=\lbrace\varphi+w_{cu}(\varphi)\mid\varphi\in C_{cu,0}^{cu}\rbrace
\end{equation*}
holds.

\subsection{Application of the Implicit Mapping Theorem}
Choose open neighborhoods $U_{c}$, $U_{u}$, and $U_{Y_{s}}$ of the origin in $C_{c}$, $C_{u}$, and $Y_{s}$, respectively, such that all subset relations
\begin{equation*}
  U_{c}+U_{Y_{s}}\subset C_{cs,0}^{cs},\quad
  U_{c}+U_{u}\subset C_{cu,0}^{cu},\quad\text{and}\quad U_{c}+U_{u}+U_{Y_{s}}\subset Y_{0}
\end{equation*}
are satisfied. Then define the map
\begin{equation*}
G:U_{c}\times U_{u}\times U_{Y_{s}}\to C_{u}\times Y_{s}
\end{equation*}
 by
\begin{equation*}
  G(\varphi^{c},\varphi^{u},\varphi^{Y_{s}}):=\begin{pmatrix}
    \varphi^{u}-P_{u}\,w_{cs}(\varphi^{c}+\varphi^{Y_{s}})\\
    \varphi^{Y_{s}}-P_{Y_{s}}w_{cu}(\varphi^{c}+\varphi^{u})
   \end{pmatrix}.
\end{equation*}
We clearly have $G(0,0,0)=(0,0)^{T}$, and as a composition of $C^{1}$-smooth maps the map $G$ is $C^{1}$-smooth as well. Now, a straightforward calculation shows that the derivative of $G$ with respect to the last two components can be represented by the matrix
\begin{equation*}
   D_{(2,3)}G(\varphi^{c},\varphi^{u},\varphi^{Y_{s}})=\begin{pmatrix}
    \id_{C_{u}}&-P_{u}Dw_{cs}(\varphi^{c}+\varphi^{Y_{s}})
    &\\
    -P_{Y_{s}}Dw_{cu}(\varphi^{c}+\varphi^{u})&\id_{Y_{s}}
  \end{pmatrix}.
\end{equation*}
In particular,
\begin{equation*}
  D_{(2,3)}G(0,0,0)=\begin{pmatrix}
    \id_{C_{u}}&0\\0&\id_{Y_{s}}
  \end{pmatrix},
\end{equation*}
 which is obviously a linear automorphism of $C_{u}\times Y_{s}$. Applying the Implicit Mapping Theorem, we obtain open neighborhoods $C_{c,0}^{c}$, $C_{u,0}^{c}$, and $C_{Y_{s},0}^{c}$ of $0$ in $U_{c}$, $U_{u}$, and $U_{Y_{s}}$, respectively, and a continuously differentiable map
 \begin{equation*}
 g:C_{c,0}^{c}\to C_{u,0}^{c}\times C_{Y_{s},0}^{c}
 \end{equation*}
  such that $g(0)=(0,0)^{T}$ and
 \begin{equation*}
   G(\varphi^{c},\varphi^{u},\varphi^{Y_{s}})=(0,0)^{T}\qquad\Longleftrightarrow\qquad g(\varphi^{c})=(\varphi^{u},\varphi^{Y_{s}})^{T}
 \end{equation*}
 for all $(\varphi^{c},\varphi^{u},\varphi^{Y_{s}})\in C_{c,0}^{c}\times C_{u,0}^{c}\times C_{Y_{s},0}^{c}$.

 \subsection{\texorpdfstring{Definition of $W_{c}$}{Definition of local center manifold}} Given $\varphi^{c}\in C^{c}_{c,0}$, set $(\varphi^{u},\varphi^{Y_{s}})^{T}:=g(\varphi^{c})$ with the map $g$ obtained in the last part. Then we have
 \begin{equation*}
 G(\varphi^{c},\varphi^{u},\varphi^{Y_{s}})=(0,0)^{T},
 \end{equation*}
  that is,
 \begin{equation*}
\left\lbrace \begin{aligned}
   0&=\varphi^{u}-P_{u}\,w_{cs}(\varphi^{c}+\varphi^{Y_{s}}),\\
   0&=\varphi^{Y_{s}}-P_{Y_{s}}w_{cu}(\varphi^{c}+\varphi^{u}).
 \end{aligned}\right.
 \end{equation*}
 By using the map $g$ and denoting by $\pi_{i}$ for $i\in\lbrace 1,2\rbrace$ the canonical projection mapping an element $(u,s)^{T}\in C_{u}\times Y_{s}$ to its $i$-th component, the last system of equations may equivalently be written as
 \begin{equation*}
   \left\lbrace
   \begin{aligned}
     (\pi_{1}\circ g)(\varphi^{c})&=P_{u}w_{cs}(\varphi^{c}+(\pi_{2}\circ g)(\varphi^{c})),\\
  (\pi_{2}\circ g)(\varphi^{c})&=P_{Y_{s}}w_{cu}(\varphi^{c}+(\pi_{1}\circ g)(\varphi^{c})).
   \end{aligned}
   \right.
 \end{equation*}
But the more important point is that on the one hand
 \begin{equation*}
   \begin{aligned}
    R\left(\varphi^{c}+\varphi^{u}+\varphi^{Y_{s}}\right)&=R\left(\varphi^{c}+\varphi^{Y_{s}}+P_{u}w_{cs}\left(
    \varphi^{c}+\varphi^{Y_{s}}\right)\right)\\
    &=(R\circ P)\left(\varphi^{c}+\varphi^{Y_{s}}+w_{cs}\left(\varphi^{c}+\varphi^{Y_{s}}\right)\right)\\
    &=(R\circ K)\left(\varphi^{c}+\varphi^{Y_{s}}+w_{cs}\left(\varphi^{c}+\varphi^{Y_{s}}\right)\right)\\
    &=\varphi^{c}+\varphi^{Y_{s}}+w_{cs}\left(\varphi^{c}+\varphi^{Y_{s}}\right)
    \in W_{cs}\cap \tilde{N}_{cs}
   \end{aligned}
 \end{equation*}
and on the other hand
 \begin{equation*}
   \begin{aligned}
    R\left(\varphi^{c}+\varphi^{u}+\varphi^{Y_{s}}\right)&=R\left(\varphi^{c}+\varphi^{u}+P_{Y_{s}}w_{cu}\left(
    \varphi^{c}+\varphi^{u}\right)\right)\\
    &=(R\circ P)\left(\varphi^{c}+\varphi^{u}+w_{cu}\left(\varphi^{c}+\varphi^{u}\right)\right)\\
    &=(R\circ K)\left(\varphi^{c}+\varphi^{u}+w_{cu}\left(\varphi^{c}+\varphi^{u}\right)\right)\\
    &=\varphi^{c}+\varphi^{u}+w_{cu}\left(\varphi^{c}+\varphi^{u}\right)\in W_{cu}\cap\tilde{N}_{cu}.
   \end{aligned}
 \end{equation*}
Thus, starting with any point $\varphi\in C_{c,0}^{c}$ we get a point, namely,
\begin{equation}\label{eq: center}
\psi=R\left(\varphi+(\pi_{1}\circ g)(\varphi)+(\pi_{2}\circ g)(\varphi)\right),
\end{equation}
 belonging to the intersection $W_{cu}\cap W_{cs}$. From now on, let $W_{c}$ denote the set of all points $\psi\in W_{cu}\cap W_{cs}$ obtained in this way.

 \subsection{\texorpdfstring{Representation of $W_{c}$ as a graph of a $C^{1}$-smooth map}{Representation of local center manifold as a graph}}
Consider the set $W_{c}$ constructed in the last step. Obviously, $W_{c}$ may alternatively be defined as the image of the map
\begin{equation*}
  C_{c,0}^{c}\ni\varphi\mapsto\varphi +(\pi_{2}\circ g)(\varphi)+w_{cs}(\varphi+(\pi_{2}\circ g)(\varphi))\in\tilde{N}_{cu}\cap \tilde{N}_{cs}
\end{equation*}
or of the map
 \begin{equation*}
   C_{c,0}^{c}\ni \varphi\mapsto \varphi+(\pi_{1}\circ g)(\varphi)+w_{cu}(\varphi+
   (\pi_{1}\circ g)(\varphi))\in \tilde{N}_{cu}\cap\tilde{N}_{cs}.
 \end{equation*}
 Therefore, $W_{c}$ particularly coincides with
the graph
 \begin{equation*}
   \lbrace \varphi+w_{c}(\varphi)\mid\varphi\in C_{c,0}^{c}\rbrace
 \end{equation*}
 of the map
 \begin{equation*}
   w_{c}:C_{c,0}^{c}\ni\varphi\mapsto (\pi_{1}\circ g)(\varphi)+w_{cu}(\varphi+(\pi_{1}\circ g)(\varphi))\in C^{1}_{su,0}
 \end{equation*}
 from the open neighborhood $C_{c,0}^{c}$ of the origin in $C_{c}$ into the open neighborhood $C_{su,0}^{1}:=C_{u,0}^{c}+C_{Y_{s},0}^{c}+C_{Z,0}^{c}$ of the origin in $C_{u}\oplus C_{s}^{1}$ where $C_{Z,0}^{c}$ denotes the open subset
 \begin{equation*}
   \left\lbrace z\in Z\mid C_{Y_{s},0}^{c}+\lbrace z\rbrace\subset C_{sz,0}^{cu},C_{u,0}^{c}+\lbrace z\rbrace\subset C_{uz,0}^{cs}\right\rbrace
 \end{equation*}
  of $Z$. In view of $g(0)=(0,0)^{T}$ and $w_{cu}(0)=0$, we also have $w_{c}(0)=0$ and so $0\in W_{c}$. Moreover, as a sum and composition of continuously differentiable maps the map $w_{c}$ is $C^{1}$-smooth as well.  We claim that
 \begin{equation*}
   Dw_{c}(0)=0.
 \end{equation*}
 In order to see this, observe that for all $\varphi\in C_{c,0}^{c}$ and all $\psi\in C_{c}$ we have
 \begin{equation*}
   Dw_{c}(\varphi)\psi=\pi_{1} Dg(\varphi)\psi+Dw_{cu}(\varphi+(\pi_{1}\circ g)(\varphi))[\id_{C_{c}}+\pi_{1} Dg(\varphi)]\psi.
 \end{equation*}
 Especially, in case $\varphi=0$
 \begin{equation*}
   Dw_{c}(0)\psi=\pi_{1}Dg(0)\psi+Dw_{cu}(0)[\id_{C_{c}}+\pi_{1} Dg(0)]\psi=\pi_{1}Dg(0)\psi
 \end{equation*}
 as $Dw_{cu}(0)=0$. Thus, for the proof of $Dw_{c}(0)=0$, it suffices to show $Dg(0)=0$. But this point is easily seen as follows. The second part implies that for all $\varphi\in C_{c,0}^{c}$ and all $\psi\in C_{c}$ we have
 \begin{equation*}
   D[\tilde{\varphi}\mapsto G\left(\tilde{\varphi},(\pi_{1}\circ g)(\tilde{\varphi}),(\pi_{2}\circ g)(\tilde{\varphi})\right)](\varphi)\psi=(0,0)^{T}\in C_{u}\times Y_{s},
 \end{equation*}
 that is, in matrix notation,
 \begin{equation*}
   \begin{aligned}
     &(0,0)^{T}\\
     =&D[\tilde{\varphi}\mapsto G(\tilde{\varphi},(\pi_{1}\circ g)(\tilde{\varphi}),(\pi_{2}\circ g)(\tilde{\varphi}))](\varphi)\,\psi\\
     =&D_{1}G(\varphi,(\pi_{1}\circ g)(\varphi),(\pi_{2}\circ g)(\varphi))\,\psi\\
     & +D_{(2,3)}[\tilde{\varphi}\mapsto G(\tilde{\varphi},(\pi_{1}\circ g)(\tilde{\varphi}),(\pi_{2}\circ g)(\tilde{\varphi}))](\varphi)\,\psi\\
     =& \begin{pmatrix}
    -P_{u}\,Dw_{cs}(\varphi+(\pi_{2}\circ g)(\varphi))\,\psi\\
    -P_{Y_{s}}\,Dw_{cu}(\varphi+(\pi_{1}\circ g)(\varphi))\,\psi
  \end{pmatrix}\\
  &+\begin{pmatrix}
    \id_{C_{u}}&-P_{u}\,Dw_{cs}(\varphi+(\pi_{2}\circ g)(\varphi))\\
    -P_{Y_{s}}Dw_{cu}(\varphi+(\pi_{1}\circ g)(\varphi))&\id_{Y_{s}}
  \end{pmatrix}Dg(\varphi)\psi.
   \end{aligned}
 \end{equation*}
 Hence, after taking into account $g(0)=(0,0)^{T}$, $Dw_{cu}(0)=0$, and $Dw_{cs}(0)=0$, for $\varphi=0$ we get
 \begin{equation*}
 \begin{aligned}
   \begin{pmatrix}0\\0\end{pmatrix}
   &=\begin{pmatrix}
     \id_{C_{u}}&0\\ 0&\id_{Y_{s}}
   \end{pmatrix} Dg(0)\psi
   &=\begin{pmatrix}\pi_{1}\left(Dg(0)\psi\right)\\
   \pi_{2}\left(Dg(0)\psi\right)\end{pmatrix}
 \end{aligned}
 \end{equation*}
for all $\psi\in C_{c}$. It follows that $Dg(0)=0\in \mathfrak{L}(C_{c},C_{u}\times Y_{s})$ and that finally proves $Dw_{c}(0)=0$ as claimed.

\subsection{Proof of Properties (C 1) -- (C 3)} Recall from the construction above that $0\in W_{c}$ and $W_{c}\subset W_{cs}\cap W_{cu}\subset X_{f}$. Moreover, $W_{c}$ is the graph of the continuously differentiable map $w_{c}:C_{c,0}^{c}\to C^{1}_{su,0}$ from the open neighborhood $C_{c,0}^{c}$ of $0$ in $C_{c}$ into the open neighborhood $C_{su,0}^{1}$ of $0$ in $C_{u}\oplus C_{s}^{1}$, and it holds that $w_{c}(0)=0$ and $Dw_{c}(0)=0$. Hence, it is clear that $W_{c}$ is not only a non-empty subset of the solution manifold $X_{f}$ but forms a $C^{1}$-smooth submanifold of $X_{f}$ with $\dim W_{c}=\dim C_{c}$.

 Next, consider the open neighborhood $N_{c}:=C_{c,0}^{c}+C^{1}_{su,0}$
of $0$ in $U_{0}$ and note that
\begin{equation*}
N_{c}\subset \tilde{N}_{cs}\cap\tilde{N}_{cu}\subset N_{cs}\cap N_{cu}.
\end{equation*}
Suppose now that $\gamma:\R\to X_{f}$ is a trajectory of $F$ with $\gamma(\R)\subset N_{c}$. We claim that for each $t\in\R$ we have $\gamma(t)\in W_{c}$. In order to see this, fix an arbitrary $T\in\R$ and set $\psi:=\gamma(T)$. Then, on the one hand, for each $t\geq 0$ we have
\begin{equation*}
  F(t,\psi)=F(t,\gamma(T))=F(t+T-T,\gamma(T))=\gamma(t+T)\in N_{c}\subset \tilde{N}_{cs}\subset N_{cs}
\end{equation*}
and therefore
\begin{equation*}
  \psi\in W_{cs}\cap\tilde{N}_{cs}
\end{equation*}
in view of property (CS 2) of the manifold $W_{cs}$. On the other hand, for the induced trajectory $\tilde{\gamma}:(-\infty,0]\ni t\mapsto \gamma(t+T)\in X_{f}$ of $F$  we have \begin{equation*}
\tilde{\gamma}(t)=\gamma(t+T)\in N_{c}\subset\tilde{N}_{cu}\subset N_{cu}
\end{equation*}
and for this reason
\begin{equation*}
\psi=\gamma(T)=\tilde{\gamma}(0)\in W_{cu}\cap \tilde{N}_{cu}
\end{equation*}
due to property (CU 1) of $W_{cu}$. Consequently,
\begin{equation*}
  \psi\in (W_{cs}\cap \tilde{N}_{cs})\cap (W_{cu}\cap\tilde{N}_{cu})
\end{equation*}
and from the local graph representations of $W_{cs}$ and $W_{cu}$ it follows that
\begin{equation*}
  \psi=P_{c}\psi+P_{Y_{s}}\psi+w_{cs}(P_{c}\psi+P_{Y_{s}}\psi)
  \quad\text{and}\quad\psi=P_{c}\psi+P_{u}\psi+w_{cu}(P_{c}\psi+P_{u}\psi).
\end{equation*}
Applying the projection operators $P_{u}$ to the first and $P_{Y_{s}}$ to the second representation of $\psi$ above, we easily infer that
\begin{equation*}
  \left\lbrace
  \begin{aligned}
    P_{u}\psi&=P_{u}\,w_{cs}(P_{c}\psi+P_{Y_{s}}\psi),\\
    P_{Y_{s}}\psi&=P_{Y_{s}}w_{cu}(P_{c}\psi+P_{u}\psi).
  \end{aligned}
  \right.
\end{equation*}
As $\psi\in N_{c}$ also implies that
\begin{equation*}
(P_{c}\psi,P_{u}\psi,P_{Y_{s}}\psi)\in C_{c,0}^{c}\times C_{u,0}^{c}\times C_{Y_{s},0}^{c}
\end{equation*}
we conclude that $G(P_{c}\psi,P_{u}\psi,P_{Y_{s}}\psi)=(0,0)^{T}$ and therefore
\begin{equation*}
  g(P_{c}\psi)=(P_{u}\psi,P_{Y_{s}}\psi)^{T}.
\end{equation*}
Hence,
\begin{equation*}
  \begin{aligned}
    \gamma(T)&=\psi\\
    &=(R\circ K)(\psi)\\
    &=R(P\psi)\\
    &=R(P_{c}\psi+P_{u}\psi+P_{Y_{s}}\psi)\\
    &=R(P_{c}\psi+(\pi_{1}\circ g)(P_{c}\psi)+(\pi_{2}\circ g)(P_{c}\psi))\in W_{c}
  \end{aligned}
\end{equation*}
due to Eq. \eqref{eq: center}. As $T\in\R$ was arbitrary chosen, the above proves our claim that $\gamma(t)\in W_{c}$ for all $t\in\R$.

Finally, we assert that, for all $\varphi\in W_{c}$ and all $\alpha>0$ such that $F(t,\varphi)$ is defined and contained in $N_{c}$ as $0\leq t<\alpha$, we also have $F(t,\varphi)\in W_{c}$ for all $0\leq t<\alpha$. Indeed, under given assumptions we clearly have
\begin{equation*}
F(t,\varphi)\in \tilde{N}_{cs}\cap \tilde{N}_{cu}\subset N_{cs}\cap N_{cu}
\end{equation*}
 for all $0\leq t<\alpha$ and therefore, by property (CS 1) of $W_{cs}$ and property (CU 1) of $W_{cu}$, $F(t,\varphi)\in W_{cs}\cap W_{cu}$ for each $0\leq t<\alpha$. Hence, it follows first that
 \begin{equation*}
   F(t,\varphi)\in (W_{cs}\cap \tilde{N}_{cs})\cap (W_{cu}\cap\tilde{N}_{cu})
 \end{equation*}
 and next, by using completely similar arguments as applied above to $\psi$, that $F(t,\varphi)\in W_{c}$ as $0\leq t<\alpha$. This finishes the proof of the assertion and completes the proof of Theorem \ref{thm: CMT}.\qed


\begin{thebibliography}{10}
\bibitem{Diekmann1995}
O. Diekmann, S. A. van Gils, S. M. Verduyn Lunel, and H.-O. Walther, {\it Delay Equations. Functional, complex, and nonlinear analysis},
\newblock Applied Mathematical Sciences 110, Springer-Verlag, New York, 1995.

\bibitem{Hartung2006}
F. Hartung, T. Krisztin, H.-O. Walther, and J. Wu, Functional differential equations
with state-dependent delay,
\newblock In: \emph{Hand. Differ. Equ.: Ordinary differential equations}, vol. \textbf{III}, Elsevier/North-Holland, Amsterdam, 2006, 435--545.

\bibitem{Kelley1967}
A. Kelley, The stable, center-stable, center,center-unstable, unstable manifolds,
\newblock \emph{J. Differ. Equations} \textbf{3} (1967), 546 -- 570.

\bibitem{Krisztin2006}
T. Krisztin, $C^{1}$-smoothness of center manifolds for differential equations with state-dependent delay,
\newblock in \emph{Nonlinear Dynamics and Evolution Equations} (eds. H. Brunner et al.), Fields Inst. Commun., \textbf{48}, Amer. Math. Soc., Providence, 2006, 213 -- 226.

\bibitem{Qesmi2009}
R. Qesmi and H.-O. Walther, Center-stable manifolds for differential equations with state-dependent delay,
\newblock \emph{Discrete Contin. Dyn. Syst.}, \textbf{23} (2009), no. 3, 1009--1033.

\bibitem{Stumpf2011}
E. Stumpf, The existence and $C^{1}$-smoothness of local center-unstable manifolds for differential equations with state-dependent delay,
\newblock \emph{Rostock. Math. Kolloq.} \textbf{66} (2011), 3--44.

\bibitem{Vanderbauwhede1989}
A. Vanderbauwhede, Centre manifolds, normal forms and elementary bifurcations,
\newblock in \emph{Dynamics Reported}, A series in dynamical systems and their applications, Vol. 2, Wiley, Chichester, 1989, 89 -- 169.

\bibitem{Walther2003}
H.-O. Walther, The solution manifold and $C^{1}$-smoothness for differential equations with state-dependent delay,
\newblock \emph{J. of Differential Equations}, \textbf{195} (2003), no. 1, 46--65.
\end{thebibliography}
\end{document}